\documentclass[10pt,twoside,reqno,draft]{amsart}

\newcommand{\vs}{\vspace}
\newcommand{\bc}{\begin{center}}
\newcommand{\ec}{\end{center}}

\usepackage{amsmath}
\usepackage{amssymb}
\usepackage{amsfonts}
\usepackage{amsthm}
\begin{document}
\setcounter{page}{1}


\vs*{1.5cm}

\title[On an equation characterizing multi-cubic
mappings]{\large On an equation characterizing multi-cubic
mappings and its stability and hyperstability}
\author[]{Abasalt Bodaghi$^{*}$ and Behrouz Shojaee$^{**}$}
\date{}
\maketitle

\vs*{-0.5cm}

\bc
{\footnotesize
$^{*}$Department of Mathematics, Garmsar Branch, Islamic Azad University, Garmsar, Iran\\
E-mail: abasalt.bodaghi@gmail.com\\
\medskip
$^{**}$Department of Mathematics, Karaj Branch, Islamic Azad University, Karaj, Iran\\
E-mail: shoujaei@kiau.ac.ir\\
}
\ec

\bigskip

{\footnotesize
\noindent
{\bf Abstract.} In this paper, we introduce $n$-variables mappings which are cubic in each variable. We show that such mappings satisfy a 
functional equation. The main purpose is to extend the applications of a fixed point method to establish the Hyers-Ulam stability for the multi-cubic mappings. As a consequence, we prove that a multi-cubic functional equation can be hyperstable.


}

\bigskip

\section{Introduction}

The study of stability problems for functional equations is related to a question of Ulam \cite{17} concerning the stability of group homomorphisms and affirmatively answered for Banach spaces by Hyers \cite{7}. Later on, various generalizations and extension of Hyers' result were ascertained by Aoki \cite{2},  Th. M. Rassias \cite{13},  J. M. Rassias \cite{ra1} and G\u{a}vru\c{t}a \cite{Gavruta} in different versions. Since then, the stability problems have been extensively investigated for a variety of functional equations and spaces. 

Let $V$ be a commutative group, $W$ be a linear space, and $n\geq 2$ be an integer. Recall from \cite{cip1} that a mapping $f: V^n\longrightarrow W$ is called {\it multi-additive} if it is additive (satisfies Cauchy's functional
equation $A(x+y)=A(x)+A(y)$) in each variable. Some facts on such mappings can be found in \cite{kuc} and many other sources. In addition, $f$ is said to be {\it multi-quadratic} if it is quadratic (satisfies quadratic functional equation $Q(x+y)+Q(x-y)=2Q(x)+2Q(y)$) in each variable  \cite{cip2}. In \cite{zha}, Zhao et al. proved that the mapping $f: V^n\longrightarrow W$ is multi-quadratic if and only if the following relation holds
\begin{align}\label{zh}
\sum_{t\in\{-1,1\}^n} f(x_1+tx_2)=2^{n}\sum_{j_1,j_2,\cdots,j_n\in \{1,2\}}f(x_{1j_1},x_{2j_2},\cdots,x_{nj_n})
\end{align}
where $x_j=(x_{1j},x_{2j},\cdots,x_{nj})\in V^n$ with $j\in\{1,2\}$. In \cite{cip1} and \cite{cip2}, Ciepli\'{n}ski studied the generalized Hyers-Ulam stability of multi-additive and multi-quadratic mappings in Banach spaces, respectively (see also \cite{zha}). 

One of the functional equations in the field of stability of functional equations is the cubic functional equation
\begin{eqnarray}\label{r}
C(x+2y)-3C(x+y)+3C(x)-C(x-y)=6C(y)
\end{eqnarray}
which is introduced by J. M. Rassias in \cite{ra1} for the first time.  It is easy to see that the mapping
$f(x)=ax^3$ satisfies (\ref{r}). Thus, every solution of the cubic functional equation (\ref{r}) is said to be a  cubic mapping. Rassias established the Ulam-Hyers stability problem for these cubic mappings. The following alternative cubic functional equation
\begin{eqnarray}\label{jk}
\mathfrak C(2x+y)+\mathfrak C(2x-y)=2\mathfrak C(x+y)+2\mathfrak C(x-y)+12\mathfrak C(x)
\end{eqnarray}
has been introduced by Jun and Kim in \cite{jk1}. They found out 
the general solution and proved the Hyers-Ulam
stability for the functional equation (\ref{jk}); for other forms of the (generalized) cubic functional equations and their stabilities on the various Banach spaces refer to \cite{bodagh}, \cite{bod1}, \cite{bmr}, \cite{jk}.

In this paper, we define multi-cubic mappings and present a characterization of such mappings. In other words, we reduce the system of $n$ equations defining the multi-cubic mappings to obtain a single equation. We also prove the generalized Hyers-Ulam stability for multi-cubic functional equations by applying the fixed point method which was introduced and used for the first time by  Brzd\c{e}k et al., in \cite{bc2}; for more applications of this approach for the satbility of multi-Cauchy-Jensen mappings in Banach spaces and 2-Banach spaces see \cite{bco2} and \cite{bc3}, respectively.

\bigskip

\section{Characterization of multi-cubic mappings}

Throughout this paper, $\mathbb N$ stands for the set of all positive integers, $\mathbb N_0:=\mathbb N\cup\{0\}, \mathbb R_+:=[0,\infty), n\in\mathbb N$. For any $l\in \mathbb N_0,m\in \mathbb N$, $t=(t_1,\cdots,t_m)\in \{-1,1\}^m$ and $x=(x_1,\cdots,x_m)\in V^m$ we write $lx:=(lx_1,\cdots,lx_m)$ and $tx:=(t_1x_1,\cdots,t_mx_m)$, where $ra$ stands, as usual, for the $r$th power of an element $a$ of the commutative group $V$.

 From now on, let $V$ and $W$ be vector spaces over the rationals, $n\in\mathbb N$ and $x_i^n=(x_{i1},x_{i2},\cdots,x_{in})\in V^n$, where $i\in\{1,2\}$. We shall denote $x_i^n$  by $x_i$ if there is no risk of ambiguity. Let $x_1,x_2\in  V^n$ and $T\in \mathbb N_0$ with $0\leq T\leq n$. Put $\mathcal M^n=\left\{(N_1,N_2, \cdots,N_n)|\,\,N_j\in\{x_{1j}\pm x_{2j},x_{1j}\}\right\}$, where $j\in\{1,\cdots,n\}$. Consider
$$\mathcal M_{T}^n :=\left\{\mathfrak N_n=(N_1,N_2, \cdots,N_n)\in \mathcal M^n|\,\,\text{Card}\{N_j:\, N_j=x_{1j}\}=T\right\}.$$
For $r\in\mathbb R$, we put $r\mathcal M_{T}^n=\left\{r\mathfrak N_n: \mathfrak N_n\in \mathcal M_{T}^n\right\}$. We say the mapping  $f:V^n\longrightarrow W$ is {\it $n$-multi-cubic} or {\it multi-cubic} if $f$ is cubic in each variable (see the equation (\ref{jk})). For such mappings, we use the following notations:
 \begin{align}\label{mtn}
 f\left(\mathcal M_{T}^n\right):=\sum_{\mathfrak N_n\in\mathcal M_{T}^n}f(\mathfrak N_n),
 \end{align}
$$f\left(\mathcal M_{T}^n,z\right):=\sum_{\mathfrak N_n\in\mathcal M_{T}^n}f(\mathfrak N_n,z)\qquad (z\in V).$$

{\bf Remark 2.1.}
It is easily verified that if the mapping $h$ satisfies the equation (\ref{jk}), then
\begin{eqnarray}\label{ab}
h(2x)=8h(x).
\end{eqnarray}
But the converse is not true. Let $(\mathcal A, \|\cdot\|)$ be a Banach algebra. Fix the vector $a_0$ in $\mathcal A$ (not necessarily unit). Define the mapping  $h:\mathcal A\longrightarrow \mathcal A$ by $h(a)=\|a\|^3a_0$ for any $a\in \mathcal A$. Clearly, for each $x\in \mathcal A$, $h(2x)=8h(x)$ while the relation (\ref{jk}) does not hold for $h$ even if we put $x=0$ and $0\neq y$. Therefore, the condition (\ref{ab}) does not imply that $h$ is a cubic mapping. 

{\bf Proposition 2.2.} {\it If the mapping $f: V^n\longrightarrow W$ is multi-cubic, then $f$ satisfies the equation 
\begin{align}\label{eqm}
\sum_{q\in\{-1,1\}^n} f(2x_1+qx_2)=\sum_{k=0}^{n}2^{n-k}12^{k}f(\mathcal M_{k}^n),
\end{align}
where $ f\left(\mathcal M_{k}^n\right)$ is defined in \eqref{mtn}. }

 \begin{proof} We prove $f$ satisfies the equation (\ref{eqm})  by induction on $n$. For $n=1$, it is trivial that $f$ satisfies the equation (\ref{jk}). If (\ref{eqm}) is valid for some
positive integer $n>1$, then,
  \begin{align*}
\nonumber \sum_{q\in\{-1,1\}^{n+1}} f(2x_1^{n+1}+qx_2^{n+1})&=2\sum_{q\in\{-1,1\}^{n}} f(2x_1^{n}+qx_2^{n},x_{1n+1}+x_{2n+1})\\
\nonumber &\,\,\,\, +2\sum_{q\in\{-1,1\}^{n}} f(2x_1^{n}+qx_2^{n},x_{1n+1}-x_{2n+1})\\
\nonumber &\,\,\,\, +12\sum_{q\in\{-1,1\}^{n}} f(2x_1^n+qx_2^n,x_{1n+1})\\
\nonumber &=2\sum_{k=0}^{n}\sum_{q\in\{-1,1\}}2^{n-k}12^{k}f(\mathcal M_{k}^n,x_{1n+1}+qx_{2n+1})\\
\nonumber &\,\,\,\, +12\sum_{k=0}^{n}2^{n-k}12^{k}f(\mathcal M_{k}^n,x_{1n+1})\\
&=\sum_{k=0}^{n+1}2^{n+1-k}12^{k}f(\mathcal M_{k}^{n+1}).
\end{align*}
 This means that (\ref{eqm}) holds for $n+1$.
\end{proof}

In the sequel,  $\left(\begin{array}{ccccc}
n\\
k\\
\end{array}\right)$ is the binomial coefficient defined for
all $n, k\in \mathbb{N}_0$ with $n\geq k$ by $n!/(k!(n-k)!)$.
   
 We say the mapping $f: V^n\longrightarrow W$ satisfies (has)  {\it the $r$-power condition} in the $j$th variable if 
$$f(z_1,\cdots,z_{j-1},2z_j,z_{j+1},\cdots, z_n)=2^rf(z_1,\cdots,z_{j-1},z_j,z_{j+1},\cdots, z_n),$$  
for all $(z_1,\cdots,z_n)\in V^n$. It follows from Remark 2.1 that the 3-power condition does not imply $f$ is cubic in the $j$th variable. Using this condition, we show that if $f$ satisfies the equation (\ref{eqm}), then it is multi-cubic as follows:      

{\bf Proposition 2.3.} {\it If the mapping $f: V^n\longrightarrow W$ satisfies the equation \emph{(\ref{eqm})} and $3$-power condition in each variable, then it is multi-cubic. }

\begin{proof} Fix $j\in\{1,\cdots,n\}$. Putting $x_{2k}=0$ for all $k\in\{1,\cdots,n\}\backslash\{j\}$ in the left side of (\ref{eqm}) and using the assumption, we get
\begin{align}\label{ea2}
\nonumber &2^{n-1}\times 2^{3(n-1)}[f\left(x_{11},\cdots,x_{1j-1},2x_{1j}+x_{2j},x_{1j+1},\cdots,x_{1n}\right)\\
\nonumber &+f\left(x_{11},\cdots,x_{1j-1},2x_{1j}-x_{2j},x_{1j+1},\cdots,x_{1n}\right)]\\
\nonumber &=2^{n-1} [f\left(2x_{11},\cdots,2x_{1j-1},2x_{1j}+x_{2j},2x_{1j+1},\cdots,2x_{1n}\right)\\
&+f\left(2x_{11},\cdots,2x_{1j-1},2x_{1j}-x_{2j},2x_{1j+1},\cdots,2x_{1n}\right)].
\end{align}
Set 
\begin{align*}
f^*(x_{1j},x_{2j}):&=f\left(x_{11},\cdots,x_{1j-1},x_{1j}+x_{2j},x_{1j+1},\cdots,x_{1n}\right)\\
&+f\left(x_{11},\cdots,x_{1j-1},x_{1j}-x_{2j},x_{1j+1},\cdots,x_{1n}\right).
\end{align*}
By the above replacements in \eqref{eqm}, it follows from (\ref{ea2}) that
\begin{align}\label{main}
\nonumber &2^{n-1}\times 2^{3(n-1)}[f\left(x_{11},\cdots,x_{1j-1},2x_{1j}+x_{2j},x_{1j+1},\cdots,x_{1n}\right)\\
 \nonumber &+f\left(x_{11},\cdots,x_{1j-1},2x_{1j}-x_{2j},x_{1j+1},\cdots,x_{1n}\right)]\\
\nonumber&=2^{n-1}\times 2^{n}f^*(x_{1j},x_{2j})\\
\nonumber&+\sum_{k=1}^{n-1}\left[\left(\begin{array}{ccccc}n-1\\k-1\\\end{array}\right)2^{2(n-k)}\times 12^{k}\right]f\left(x_{11},\cdots,x_{1n}\right)\\
\nonumber&+\sum_{k=1}^{n-1}\left[\left(\begin{array}{ccccc}n-1\\k\\\end{array}\right)2^{2(n-k)-1}\times 12^{k}\right]f^*(x_{1j},x_{2j})\\
\nonumber&+12^{n}f\left(x_{11},\cdots,x_{1n}\right)\\
\nonumber&=\left[2^{2n-1}+\sum_{k=1}^{n-1}\left(\begin{array}{ccccc}n-1\\k\\\end{array}\right)2^{2(n-k)-1}\times 12^{k}\right]f^*(x_{1j},x_{2j})\\
&+\left[12^{n}+\sum_{k=1}^{n-1}\left(\begin{array}{ccccc}n-1\\k-1\\\end{array}\right)2^{2(n-k)}\times 12^{k}\right]f\left(x_{11},\cdots,x_{1n}\right).
\end{align}
On the other hand, we have
\begin{align}\label{w2}
\nonumber 2^{2n-1}+\sum_{k=1}^{n-1}\left(\begin{array}{ccccc}n-1\\k\\\end{array}\right)2^{2(n-k)-1}\times 12^k&= 2^{2n-1}\left(1+\sum_{k=1}^{n-1}\left(\begin{array}{ccccc}n-1\\k\\\end{array}\right)3^k\right)\\
&=2^{2n-1}\left(1+3\right)^{n-1}=2^{4n-3}.
 \end{align}
 In addition,
  \begin{align}\label{w1}
\nonumber 12^n+\sum_{k=1}^{n-1}\left(\begin{array}{ccccc}n-1\\k-1\\\end{array}\right)2^{2(n-k)}\times 12^k&=12^n+\sum_{k=1}^{n-1}\left(\begin{array}{ccccc}n-1\\k-1\\\end{array}\right)2^{2(n-k)}\times 2^{2k}\times 3^{k}\\
\nonumber &=12^n+3\times 2^{2n}\sum_{k=0}^{n-2}\left(\begin{array}{ccccc}n-1\\k-1\\\end{array}\right) 3^{k}\\
\nonumber &=12^n+3\times 2^{2n}\left(\sum_{k=0}^{n-1}\left[\left(\begin{array}{ccccc}n-1\\k-1\\\end{array}\right) 3^{k}\right]-3^{n-1}\right)\\
\nonumber &=12^n+3\times 2^{2n}\left((1+3)^{n-1}-3^{n-1}\right)\\
\nonumber &=12^n+3\times 2^{2n}\left(2^{2(n-1)}-3^{n-1}\right)\\
&=12\times 2^{4(n-1)}.
 \end{align}
The relations (\ref{main}), (\ref{w2}) and (\ref{w1}) imply that
\begin{align*}
&f\left(x_{11},\cdots,x_{1j-1},2x_{1j}+x_{2j},x_{1j+1},\cdots,x_{1n}\right)\\
 &+f\left(x_{11},\cdots,x_{1j-1},2x_{1j}-x_{2j},x_{1j+1},\cdots,x_{1n}\right)\\
&=2f^*(x_{1j},x_{2j})+12f\left(x_{11},\cdots,x_{1n}\right)
\end{align*}
This means that $f$ is cubic in the $j$th  variable. Since $j$ is arbitrary, we obtain the desired result. 
 \end{proof}
 

\section{Stability Results for (\ref{eqm})}
In this section, we prove the generalized Hyers-Ulam stability of equation (\ref{eqm}) by a fixed point result (Theorem 3.1) in Banach spaces. Throughout, for two sets $X$ and $Y$, the set of all mappings from $X$ to $Y$  is denoted by $Y^X$. We introduce the upcoming three hypotheses:

\begin{enumerate}
\item[(A1)] {$Y$ is a Banach space, $\mathcal S$ is a nonempty set, $j\in\mathbb N$, $g_1,\cdots,g_j:\mathcal S\longrightarrow \mathcal S$ and $L_1,\cdots,L_j:\mathcal S\longrightarrow \mathbb R_+$,}
\item[(A2)] {$\mathcal T:Y^\mathcal S\longrightarrow Y^\mathcal S$ is an operator satisfying the inequality
$$\left\|\mathcal T\lambda(x)-\mathcal T\mu(x)\right\|\leq \sum_{i=1}^jL_i(x)\left\|\lambda(g_i(x))-\mu(g_i(x))\right\|,\quad \lambda,\mu\in Y^\mathcal S, x\in \mathcal S,$$}
\item[(A3)] {$\Lambda:\mathbb R_+^\mathcal S\longrightarrow \mathbb R_+^\mathcal S$ is an operator defined through
$$\Lambda\delta(x):=\sum_{i=1}^jL_i(x)\delta(g_i(x))\qquad \delta\in \mathbb R_+^\mathcal S, x\in \mathcal S.$$}
\end{enumerate}

Here, we highlight the following theorem which is a fundamental result in fixed point theory \cite[Theorem 1]{bc2}. This result plays a key tool to obtain our objecive in this paper.

{\bf Theorem 3.1.}
{\it Let hypotheses \emph{(A1)-(A3)} hold and the function $\theta:\mathcal S\longrightarrow \mathbb R_+$ and the mapping  $\phi:\mathcal S\longrightarrow Y$ fulfill the following two conditions:
$$\|\mathcal T\phi(x)-\phi(x)\|\leq \theta(x),\quad \theta^*(x):=\sum_{l=0}^\infty\Lambda^l\theta(x)<\infty\qquad (x\in \mathcal S).$$
Then, there exists a unique fixed point $\psi$ of $\mathcal T$ such that
$$\|\phi(x)-\psi(x)\|\leq \theta^*(x)\qquad (x\in \mathcal S).$$
Moreover, $\psi(x)=\lim_{l\rightarrow\infty}\mathcal T^l\phi(x)$ for all $x\in \mathcal S$.}

Here and subsequently, for the mapping $f:V^n \longrightarrow W$, we consider the difference operator $\mathfrak Df:V^n\times V^n \longrightarrow W$ by
\begin{align*}
\mathfrak Df(x_1,x_2)&:=\sum_{q\in\{-1,1\}^n} f(2x_1+qx_2)-\sum_{k=0}^{n}2^{n-k}12^{k}f\left(\mathcal M_{k}^n\right),
\end{align*}
where $f\left(\mathcal M_{k}^n\right)$ is defined in \eqref{mtn}.  With this notation, we have the next stability result for the functional equation (\ref{eqm}). 
 
{\bf Theorem 3.2.} {\it Let $\beta\in \{-1,1\}$, let $V$ be a linear space and $W$ be a Banach space. Suppose that $\phi:V^n\times V^n \longrightarrow \mathbb R_+$  is a mapping satisfying
\begin{align}\label{a00}
\lim_{l\rightarrow\infty}\left(\frac{1}{2^{3n\beta}}\right)^l \phi (2^{\beta l}x_1,2^{\beta l}x_2)=0
\end{align}
for all $x_1,x_2\in V^n$ and
\begin{align}\label{a2}
\Phi(x)=\frac{1}{2^{3n\frac{\beta+1}{2}+n}}\sum_{l=0}^{\infty}\left(\frac{1}{2^{3n\beta}}\right)^l\phi\left(2^{\beta l+\frac{\beta-1}{2}}x,0\right)< \infty
\end{align}
for all $x\in V^n$. Assume also $f:V^n \longrightarrow W$ is a mapping satisfying the inequality 
 \begin{align}\label{a1}
\|\mathfrak Df(x_1,x_2)\|_{Y}\leqslant \phi (x_1,x_2)
\end{align}
for all $x_1,x_2\in V^n$. Then, there exists a unique multi-cubic mapping $\mathcal C:V^n \longrightarrow W$ such that
\begin{equation}\label{a3}
\|f(x)-\mathcal C(x) \| \leq \Phi(x)
\end{equation} 
for all $x\in V^n$.}

\begin{proof}
Putting $x=x_1$ and $x_{2}=0$ in (\ref{a1}), we have
 \begin{align}\label{a4}
\left\|2^nf(2x)-\left(\sum_{k=0}^{n}\left(\begin{array}{ccccc}n\\k\\\end{array}\right)2^{2(n-k)}\times 12^k\right)f(x)\right\|\leq\phi(x,0)
\end{align}
for all $x\in V^n$. By an easy computation, we have
\begin{align}\label{a9p}
&\sum_{k=0}^{n}\left(\begin{array}{ccccc}n\\k\\\end{array}\right)2^{2(n-k)}\times 12^k=2^{2n}\sum_{k=0}^{n}\left(\begin{array}{ccccc}n\\k\\\end{array}\right)3^k=2^{2n}(1+3)^n=2^{4n}.
 \end{align}
It follows from (\ref{a4}) and (\ref{a9p}) that 
 \begin{align}\label{aaa}
\left\|f(2x)-2^{3n}f(x)\right\|\leq\frac{1}{2^n}\phi(x,0)
\end{align}
for all $x\in V^n$. Set  $$\xi(x):=\frac{1}{2^{3n\frac{\beta+1}{2}+n}}\phi\left(2^{\frac{\beta-1}{2}}x,0\right),\, \text{and}\, \mathcal T\xi(x):=\frac{1}{2^{3n\beta}}\xi(2^{\beta}x)\qquad (\xi\in W^{V^n}).$$
Then, the relation (\ref{aaa}) can be modified as 
\begin{align}\label{a10}
\left\|f(x)-\mathcal Tf(x)\right\|\leq\xi(x)\qquad (x\in V^n).
\end{align}
 Define $\Lambda\eta(x):=\frac{1}{2^{3n\beta}}\eta(2^{\beta}x)$ for all $\eta\in \mathbb R_+^{V^n}, x\in V^n$. We now see that $\Lambda$ has the form described in (A3) with $\mathcal S=V^n$, $g_1(x)=2^{\beta}x$ and $L_1(x)=\frac{1}{2^{3n\beta}}$ for all $x\in V^n$. Furthermore, for each $\lambda,\mu\in W^{V^n}$ and $x\in V^n$, we get
\begin{align*}
\left\|\mathcal T\lambda(x)-\mathcal T\mu(x)\right\|=\left\|\frac{1}{2^{3n\beta}}\left[\lambda(2^{\beta}x)-\mu(2^{\beta}x)\right]\right\|\leq L_1(x)\left\|\lambda(g_1(x))-\mu(g_1(x))\right\|.
\end{align*}
The above relation shows that the hypotheis (A2) holds. By induction on $l$, one can check that for any $l\in \mathbb N_0$ and $x\in V^n$, we have
\begin{align}\label{a11}
\Lambda^l\xi(x):=\left(\frac{1}{2^{3n\beta}}\right)^l\xi(2^{\beta l}x)=\frac{1}{2^{3n\frac{\beta+1}{2}+n}}\left(\frac{1}{2^{3n\beta}}\right)^l\phi\left(2^{\beta l+\frac{\beta-1}{2}}x,0\right)
\end{align}
for all $x\in V^n$. The relations (\ref{a2}) and (\ref{a11}) necessitate that all assumptions of Theorem 3.1 are satisfied. Hence, there  exsits a unique mapping $\mathcal C:V^n \longrightarrow W$ such that
$$\mathcal C(x)=\lim_{l\rightarrow\infty}(\mathcal T^lf)(x)=\frac{1}{2^{3n\beta}}\mathcal C(2^{\beta}x)\qquad(x\in V^n),$$
and (\ref{a3}) holds. We shall to show that
\begin{align}\label{a12}
\|\mathfrak D(\mathcal T^lf)(x_1,x_2)\|\leq \left(\frac{1}{2^{3n\beta}}\right)^l \phi (2^{\beta l}x_1,2^{\beta l}x_2)
\end{align}
for all $x_1,x_2\in V^n$ and $l\in\mathbb N_0$. We argue by induction on $l$. The inequality (\ref{a12}) is valid for $l=0$ by (\ref{a1}). Assume that (\ref{a12}) is true for an $l\in\mathbb N_0$. Then
\begin{align}\label{a13}
\nonumber&\|\mathfrak D(\mathcal T^{l+1}f)(x_1,x_2)\|\\
\nonumber&=\left\|\sum_{q\in\{-1,1\}^n} (\mathcal T^{l+1}f)(2x_1+qx_2)-\sum_{k=0}^{n}2^{n-k}12^{k}(\mathcal T^{l+1}f)(\mathcal M_k^n)\right\|\\
\nonumber&=\frac{1}{2^{3n\beta}}\left\|\sum_{q\in\{-1,1\}^n} (\mathcal T^{l}f)(2^{\beta}(2x_1+qx_2))-\sum_{k=0}^{n}2^{n-k}12^{k}(\mathcal T^{l}f)(2^{\beta}\mathcal M_k^n)\right\|\\
&=\frac{1}{2^{3n\beta}}\left\|\mathfrak D(\mathcal T^lf)(2^{\beta}x_1,2^{\beta}x_2)\right\|\leq \left(\frac{1}{2^{3n\beta}}\right)^{l+1} \phi (2^{\beta (l+1)}x_1,2^{\beta (l+1)}x_2)
\end{align}
for all $x_1,x_2\in V^n$. Letting $l\rightarrow\infty$ in (\ref{a12}) and applying (\ref{a00}), we arrive at $\mathfrak D\mathcal C(x_1,x_2)=0$ for all $x_1,x_2\in V^n$. This means that the mapping $\mathcal C$ satisfies (\ref{eqm}).  Finally, assume that $\mathcal C':V^n \longrightarrow W$ is another multi-cubic mapping satisfying the equation (\ref{eqm}) and inequality (\ref{a3}), and fix $x\in V^n$, $j\in \mathbb N$. Then
\begin{align*}
&\|\mathcal C(x)-\mathcal C'(x)\|\\
&=\left\|\frac{1}{2^{3n\beta j}}\mathcal C(2^{\beta j}x)-\frac{1}{2^{3n\beta j}}\mathcal C'(2^{\beta j}x)\right\|\\
&\leq\frac{1}{2^{3n\beta j}}(\|\mathcal C(2^{\beta j}x)-f(2^{\beta j}x)\|+\|\mathcal C'(2^{\beta j}x)-f(2^{\beta j}x)\|)\\
 &\leq\frac{1}{2^{3n\beta j}}2\Phi(2^{\beta j}x)\\
&\leq 2 \frac{1}{2^{3n\frac{\beta+1}{2}+n}}\sum_{l=j}^{\infty}\left(\frac{1}{2^{3n\beta}}\right)^l\phi\left(2^{\beta l+\frac{\beta-1}{2}}x,0\right).
\end{align*}
Consequently, letting $j\rightarrow\infty$ and using the fact that series (\ref{a2}) is convergent for all $x\in V^n$, we obtain $\mathcal C(x)=\mathcal C'(x)$ for all $x\in V^n$, which finishes the proof. 
\end{proof}

Let $A$ be a nonempty set, $(X,d)$ a metric space, $\psi\in \mathbb R_{+}^{A^n}$, and $\mathcal F_1, \mathcal F_2$
operators mapping a nonempty set $D\subset X^A$ into $X^{A^n}$. We say that operator equation 
\begin{align}\label{hyper}
\mathcal F_1\varphi(a_1,\cdots,a_n)=\mathcal F_2\varphi(a_1,\cdots,a_n)
\end{align}
is $\psi$-hyperstable provided every $\varphi_0\in D$ satisfying inequality 
$$d(\mathcal F_1\varphi_0(a_1,\cdots,a_n),\mathcal F_2\varphi_0(a_1,\cdots,a_n))\leq \psi(a_1,\cdots,a_n), \qquad a_1,\cdots,a_n\in A,$$
fulfils (\ref{hyper}); this definition is introduced in \cite{brc}. In other words, a functional equation $\mathcal F$ is {\it hyperstable} if any mapping $f$ satisfying the equation $\mathcal F$ approximately is a true solution of $\mathcal F$.

Under some conditions the functional equation (\ref{eqm}) can be hyperstable as follows. 

{\bf Corollary 3.4.} {\it Let $\delta>0$. Suppose that $p_{ij}>0$ for $i\in\{1,2\}$ and $j\in\{1,\cdots,n\}$ fulfill $\sum_{i=1}^{2}\sum_{j=1}^{n}p_{ij}\neq 3n$. Let $V$ be a normed space and let $W$ be a Banach space. If $f:V^n \longrightarrow W$ is a mapping satisfying the inequality 
 \begin{align*}
\|\mathfrak Df(x_1,x_2)\|\leq \prod_{i=1}^{2}\prod_{j=1}^{n}\|x_{ij}\|^{p_{ij}}\delta
\end{align*}
for all $x_1,x_2\in V^n$, then $f$ is multi-cubic.}

In the following corollary, we show that the functional equation (\ref{eqm}) is stable. Since the proof is routine, we include it without proof.
 
{\bf Corollary 3.5.} 
{\it Let $\delta>0$ and $\alpha\in \mathbb R$ with $\alpha\neq 3n$. Let also $V$ be a normed space and let $W$ be a Banach space. If $f:V^n \longrightarrow W$ is a mapping satisfying the inequality 
 \begin{align*}
\|\mathfrak Df(x_1,x_2)\|\leq \sum_{i=1}^{2}\sum_{j=1}^{n}\|x_{ij}\|^{\alpha}\delta
\end{align*}
for all $x_1,x_2\in V^n$, then there exists a unique multi-cubic mapping $\mathcal C:V^n \longrightarrow W$ such that

$$\|f(x)-\mathcal C(x)\|
 \leq \begin{cases}
\frac{\delta}{2^{4n}-2^{\alpha+n}}\sum_{j=1}^{n}\|x_{1j}\|^{\alpha}\,\,\hspace{2cm}  \alpha< 3n\\\\

\frac{2^{\alpha}}{2^{\alpha+n}-2^{4n}}\delta\sum_{j=1}^{n}\|x_{1j}\|^{\alpha}\,\,\hspace{1.8cm}  \alpha>3n
\end{cases}
$$
for all $x\in V^n$. }


\section*{Acknowledgements}
The authors sincerely appreciate the anonymous reviewer for her/his
careful reading, constructive comments and fruitful suggestions
to improve the paper.

\bigskip

\end{document}